\documentstyle[amscd,amssymb,verbatim,12pt,diagrams]{amsart}

\newcommand{\Tor}{\operatorname{Tor}}

\newcommand{\OO}{{\cal O}}

\newcommand{\codim}{\operatorname{codim}}

\renewcommand{\P}{{\Bbb P}}

\numberwithin{equation}{section}

\newtheorem{thm}{Theorem}[section]
\newtheorem{prop}[thm]{Proposition}
\newtheorem{prop-def}[thm]{Proposition-Definition}
\newtheorem{lem}[thm]{Lemma}
\newtheorem{cor}[thm]{Corollary}
\newenvironment{rem}{\vspace{3mm}\noindent
{\bf Remark.}}{\vspace{3mm}}
\newenvironment{defi}{\vspace{3mm}\noindent
{\bf Definition.}}{\vspace{3mm}}

\newcommand{\Pf}{\noindent {\it Proof}}

\newcommand{\JJ}{{\cal J}}

\newcommand{\LL}{{\cal L}}

\newcommand{\Z}{{\Bbb Z}}

\newcommand{\La}{\Lambda}

\newcommand{\wt}{\widetilde}

\newcommand{\sub}{\subset}
\newcommand{\ed}{\qed\vspace{3mm}}

\renewcommand{\span}{\operatorname{span}}

\title{Koszul configurations of points in projective spaces}
\author{A. Polishchuk}
 \thanks{Supported in part by NSF grant}

\begin{document}
\begin{abstract} We prove a new criterion for the homogeneous coordinate ring of a finite set
of points in $\P^n$ to be Koszul. Like the well known criterion due to Kempf \cite{K} it involves
only incidence conditions on linear spans of subsets of a given set. We also give a sufficient condition
for the Koszul property to be preserved when passing to a subset of a finite set of points in
$\P^n$.
\end{abstract}
\maketitle

\bigskip

\centerline{\sc Introduction}

\medskip

This paper is an attempt to understand better which finite sets of points
in $\P^n$ are {\it Koszul}, i.e., have Koszul homogeneous coordinate algebra.
Recall that Koszulness is a remarkable refininement
of the property of a graded algebra to be defined by degree one generators
and quadratic relations (see \cite{Pr}, \cite{BF}, \cite{BGS}, \cite{PP}). 
It is well known that a homogeneous coordinate algebra of a non-degenerate projectively normal 
variety $X\sub\P^n$ is Koszul iff the same is true for a hyperplane section of $X$. 
Because of this one can often reduce the question of Koszulness of a homogeneous coordinate
algebra of $X$ to that of a finite set of points. Thus, it is important to know which
configurations of points are Koszul. Theorem of Kempf \cite{K}
states that every set of $\le 2n$ points in general linear
position in $\P^n$ is Koszul. One can ask whether there are any other Koszulness criteria involving only
incidence conditions on linear spans of subsets of our set.
Kempf's theorem is the optimal criterion of this kind for sets of points in general linear position: 
if such a set contains more than $2n$ elements then it is not necessarily Koszul (see \cite{CRV}, Sec.~4, for
an example of a non-Koszul set of $9$ points in general linear position in $\P^4$).
Our main result gives a Koszulness criterion for sets of points that are not necessarily in general linear position.

Let $S\sub\P^n$ be a finite set of distinct points. We denote by $\LL(S)$ the lattice of all subsets of $S$
obtained as intersections of $S$ with linear subspaces in $\P^n$ (for example, the empty set is
an element of $\LL(S)$).

\begin{thm}\label{mainthm}
Assume there exists a subset $\La\sub\LL(S)$ such that
$|\La|>1$ and the following property holds:

\noindent
$(\ast)$ for every $T\in\La$ such that $T\neq S$ there exists $T'\in\La$ such that $T'\neq\emptyset$,  
$T'$ is linearly independent, $T\cap T'=\emptyset$, and $T\cup T'$ is an element of $\La$.

Then $S$ is Koszul. Moreover, in this situation a stronger property holds: the homogeneous coordinate
ring of $S$ admits a {\it Koszul filtration} (see \cite{CTV}).
\end{thm}

Note that Kempf's criterion of Koszulness is an easy corollary of Theorem \ref{mainthm}. 
In fact, we obtain the following generalization
of Kempf's theorem that can also be proved by methods of \cite{CTV}. 
For a finite subset $S\in\P^n$ let us denote by $\span(S)$ the linear span
of $S$ (it is a linear subspace of $\P^n$).

\begin{cor}\label{maincor} Assume that $S=S_1\sqcup S_2$ with $S_1$ and $S_2$ linearly
independent and such that 
$$\span (S_1)\cap S_2=\span(S_2)\cap S_1=\emptyset.$$ 
Then $S$ is Koszul.
\end{cor}

Indeed, in this case one can take $\La=\{S_1,S_2,S\}$. Other examples of configurations satisfying
condition $(\ast)$ will be considered in section \ref{main-sec}. 

For a subspace $V\sub H^0(\P^n,\OO(d))$ and a set of points $S\sub\P^n$ we say that
$S$ imposes independent conditions on $V$ if the evaluation map
$V\to H^0(S,\OO(d)|_S)$ is surjective. In the case $V=H^0(\P^n,\OO(2))$ we say that $S$ imposes
independent conditions on quadrics.
In the situation of Theorem \ref{mainthm} we also check that $S$ imposes independent
conditions on quadrics (see Lemma \ref{tech-lem}(a)). 
Koszulness of configurations of points that impose dependent
conditions on quadrics seems to be much more difficult to grasp. One result in this
direction is Thm. 5.3 of \cite{CRV} that involves some assumptions of the different nature than
the ones we have considered. 
The following theorem also deals with sets that may impose dependent conditions on quadrics.
It gives a criterion of Koszulness for a subset of a Koszul set of points.

\begin{thm}\label{descent-points-thm} 
Let $S$ be a Koszul set of distinct points in $\P^n$ and let $S'\sub S$ be a subset.
Assume that $S\setminus S'$ imposes independent conditions on quadrics passing through
$S'$. Then $S'$ is also Koszul.
\end{thm}

\begin{cor}\label{descent-cor}
Let $S$ be a Koszul set of points in $\P^n$ imposing independent conditions
on quadrics. Then every subset of $S$ is also Koszul.
\end{cor}

Theorem \ref{descent-points-thm} seems to indicate that in the case of a set $S$ imposing dependent
conditions on quadrics one can try to search for Koszulness criteria
involving incidence conditions on spans of various subsets of $S$ along with conditions of the kind
``$S_1$ imposes independent conditions on quadrics passing through $S_2$" for $S_1\sub S$ and
$S_2\sub S$. 

Another interesting question is for which projective configurations of points the corresponding homogeneous
ideal admits a quadratic Gr\"obner basis. It is known that this is the case for configurations described
in Corollary \ref{maincor} (see \cite{CRV}).

The paper is organized as follows. In section \ref{gen-sec} we recall two general methods
of proving Koszulness of algebras: the first is based on the principle that Koszulness is inherited
under certain types of homomorphisms and the second uses Koszul filtrations. 
Theorem \ref{descent-points-thm} is an easy application of the first method. 
In section \ref{main-sec} we prove Theorem \ref{mainthm} constructing a Koszul filtration on the
corresponding algebra.

\noindent
{\it Conventions}. Throughout this paper we work over a fixed ground field $k$.
By {\it points} in $\P^n$ we mean $k$-points.
By a {\it graded algebra} we always mean an associative $\Z$-graded $k$-algebra of the form
$A=\oplus_{n\ge 0}A_n$, where $A_0=k$.

\noindent
{\it Acknowledgment}. I am grateful to Pavel Etingof and Eric Rains for useful discussions.
Also, I thank the referee for providing a simpler proof of Theorem \ref{descent-points-thm} and other
useful remarks.

\section{General criteria of Koszulness}\label{gen-sec}

In this section for the convenience of the reader we recall some well-known methods of proving Koszulness. 
For purposes of this work we only need the case of commutative algebras. However, as the reader may
notice, this assumption is not really used (one just has to work with left ideals and left modules everywhere). 

\medskip

\noindent
{\bf Homomorphisms and Koszulness.}

\medskip

Assume we are given a surjective homomorphism $f:A\to B$ of graded algebras. 
Under certain conditions on $f$ one can establish a relation between Koszulness of
$A$ and $B$.
The well known example is the case when  $B$ admits a linear
free resolution as a (left) graded $A$-module, i.e., a resolution of the form
$$\ldots\to V_3\otimes A(-3)\to V_2\otimes A(-2)\to V_1\otimes A(-1)\to A,$$ 
where $V_i$ are vector spaces. 

\begin{prop}\label{hom-prop} 
Assume that $B$ is Koszul and admits a linear free resolution as a left $A$-module.
Then $A$ is also Koszul.
\end{prop}

\Pf .
Let us consider the standard spectral sequence 
\begin{equation}\label{standard-sp-seq}
E^2_{p,q}=\Tor_p^B(\Tor_q^A(k,B),k)\implies \Tor_{p+q}^A(k,k).
\end{equation}
The existence of a linear resolution for $B$ over $A$ implies that $\Tor_q^A(k,B)$
is always concentrated in internal degree $q$. By Koszulness of $B$ it follows
that $E_2^{pq}$ is concentrated in internal degree $p+q$. Therefore, the spectral
sequence implies that $\Tor_i^A(k,k)$ has internal degree $i$.
\ed

The converse to the above proposition is also true. In fact, one has the following stronger
result. Recall that $\Tor_{i,j}^A(k,M)$ denotes the component of internal degree $j$ in
$\Tor_i^A(k,M)$.

\begin{thm}\label{hom-thm} (\cite{Pos}, Thm.~5)
Let $A\to B$ be a homomorphism of graded algebras 
such that $\Tor_{i,j}^A(k,B)=0$ for $j>i+1$.
Assume that $A$ is Koszul. Then $B$ is also Koszul.
\end{thm}

\Pf . Let us prove by induction in $n$ that 
$\Tor_n^B(k,k)$ has internal degree $n$. 
This is true for $n=0$. Assume that the assertion is true for all 
$n'<n$ and let us show that it is true for $n$.
Consider again the spectral sequence \eqref{standard-sp-seq}.
The assumption of the theorem implies that $\Tor_i^A(k,B)$ is an extension of
trivial $B$-modules concentrated in degree $i$ and $i+1$. 
Hence, from the induction assumption we get that for $p<n$ the term $E^2_{p,q}$ is concentrated in degree 
$p+q$ and $p+q+1$.
Now let us consider the following piece of the term $E^N$ for $N\ge 2$:
$$0=E^N_{n+N,-N+1}\stackrel{d_N}{\to} E^N_{n,0}\stackrel{d_N}{\to} E^N_{n-N,N-1}.$$
Since $E^N_{n-N,N-1}$ is a subquotient of $E^2_{n-N,N-1}$, 
it is concentrated in degree $n-1$ and $n$. 
On the other hand, the limit term is concentrated in degree $n$ (since $A$ is Koszul), 
so it follows that $E^2_{n,0}=\Tor_n^B(\Tor_0^A(k,B),k)$ is concentrated in degree $n-1$ and $n$.

Next, we observe that $\Tor_0^A(k,B)$ is concentrated in degree $0$ and
$1$, so we have an exact triple of right $B$-modules
$$0\to\Tor_{0,1}^A(k,B)\otimes k(-1)\to\Tor_0^A(k,B)\to k\to 0$$
giving rise to a long exact sequence
$$\ldots\to \Tor_n^B(\Tor_0^A(k,B),k)\to \Tor_n^B(k,k)
\to \Tor_{0,1}^A(k,B)\otimes\Tor_{n-1}^B(k(-1),k)\to\ldots$$
As we have shown above the first term lives in degree $\le n$. Also,
by the induction assumption the space $\Tor_{n-1}^B(k(-1),k)$ has internal degree $n$.
Hence, $\Tor_{n,j}^B(k,k)=0$ for $j>n$. 
\ed

For example, the quotient of a Koszul algebra by a central element of degree $1$ or $2$ that is not a zero divisor, is again a Koszul algebra. Also,
the quotient of a Koszul algebra $A$ by an ideal $I\sub A$ such that $I(2)$ admits a linear free resolution, is Koszul (cf. \cite{BF}, Theorems 4(e) and 7(a); \cite{ST}, Theorems 1.2 and 1.5).

We are going to deduce Theorem \ref{descent-points-thm} from Theorem \ref{hom-thm}.
Let $S_1$ and $S_2$ be finite subsets of points in $\P^n$, and let $V_1\sub H^0(\P^n,\OO(2))$
(resp., $V_2\sub H^0(\P^n,\OO(2))$) be the space of quadrics vanishing on $S_1$ (resp., $S_2$). 
We say that $S_1$ and $S_2$ are {\it relatively $2$-independent} if $V_1+V_2=H^0(\P^n,\OO(2))$.
For example, one can easily check that if $S_1$ imposes independent conditions on quadrics
passing through $S_2$ then $S_1$ and $S_2$ are relatively $2$-independent.
Thus, Theorem \ref{descent-points-thm} is an immediate corollary of the following result.

\begin{prop} 
Let $S$ be a Koszul set of distinct points in $\P^n$ and let $S'\sub S$ be a subset.
Assume that $S'$ and $S\setminus S'$ are relatively $2$-independent. Then
$S'$ is also Koszul.
\end{prop}

\Pf . Set $S''=S\setminus S'$. Let $I$, $I'$ and $I''$ be the homogeneous ideals in 
$R=k[x_0,\ldots,x_n]$ corresponding to $S$, $S'$ and $S''$, respectively, and let
$A=R/I$, $A'=R/I'$ and $A''=R/I''$ be the corresponding homogeneous coordinate algebras.
By definition $I=I'\cap I''$, so we have an exact sequence of $A$-modules
$$0\to A\to A'\oplus A''\to M\to 0$$
where $M=R/(I'+I'')$. 
The condition that $S'$ and $S''$ are relatively $2$-independent implies that $M$ is concentrated in
degree $0$ and $1$. Therefore, $\Tor_i^A(k,M)$ is concentrated in internal degree $i$ and $i+1$ (by Koszulness
of $A$). Now the long exact sequence of $\Tor$'s shows that $\Tor_i^A(k,A')$ is also concentrated in degree $i$
and $i+1$. Applying Theorem \ref{hom-thm} we deduce that $A'$ is Koszul. 
\ed

\begin{rem} In the case when $|S'|=|S|-1$ the converse of the above proposition (or Theorem \ref{descent-points-thm})
is also true: if $S$ and $S'$ are Koszul then there exists a quadric passing through $S'$
but not passing through the point $S\setminus S'$. In fact, it suffices to assume for this that the coordinate
rings $A_S$ and $A_{S'}$ are quadratic. 
\end{rem}

\medskip

\noindent
{\bf Koszul families of ideals.}

\medskip

Now let us describe another powerful way to derive Koszulness.

\begin{defi} Let $\JJ$ be a finite family of homogeneous ideals in a graded algebra $A$.
We say that $\JJ$ is {\it Koszul} if for every $J\in \JJ$ such that $J\neq 0$,
there exists $J_1,J_2\in\JJ$ such that
$J_1\sub J$, $J_1\neq J$, and $J/J_1\simeq A/J_2(-1)$ as graded $A$-modules.
\end{defi}

Note that a nonempty Koszul family of ideals necessarily contains the zero ideal. 
The above definition is closely related to the notion of {\it Koszul filtration} considered in
\cite{CTV} and \cite{Piont}. In fact, a Koszul family gives a Koszul filtration provided 
$A_+$ is contained in $\JJ$. The following result is the main reason for introducing such
families: in the case of a Koszul filtration it leads to Koszulness of $A$.

\begin{prop}\label{Koszul-fam-prop} 
Let $\JJ$ be a Koszul family of ideals in $A$. Then for every $J\in A$ the
$A$-module $A/J$ admits a linear free resolution.
\end{prop}

\Pf . Let us prove by induction in $n$ that for every $J\in\JJ$
one has $\Tor_{i,j}^A(k,A/J)=0$ for $i\le n$ and $j>i$. For $n=0$ this is
clear. Now let us assume that $n>0$ and the assertion is true for $n-1$. 
If $J=0$ then the assertion is trivial, so we can assume that $J\neq 0$.
Then there exists $J_1,J_2\in\JJ$ such that $J_1\sub J$, $J_1\neq J$, and an isomorphism
$J/J_1\simeq A/J_2(-1)$. Therefore, we have an exact
sequence of $A$-modules
$$0\to A/J_2(-1)\to A/J_1\to A/J\to 0$$
that gives a long exact sequence
$$\ldots\to \Tor_{n}^A(k, A/J_1)\to\Tor_{n}^A(k,A/J)\to
\Tor_{n-1}^A(k,A/J_2(-1))\to\ldots$$
By induction assumption the space $\Tor_{n-1}^A(k,A/J_2(-1))$ is concentrated in
internal degree $n$. Hence, the statement reduces to the similar statement for
$J_1$. Since a strictly decreasing chain of elements of $\JJ$ has to terminate (at the zero ideal) this
completes the proof of the induction step.
\ed

\begin{cor}\label{Koszul-fam-cor} 
Assume that $A$ admits a Koszul family of ideals $\JJ$ such that for some
$J\in \JJ$ the algebra $A/J$ is Koszul. Then the algebra $A$ is also Koszul.
\end{cor}

\Pf . This follows from Propositions \ref{hom-prop} and \ref{Koszul-fam-prop}.
\ed

In particular, if a Koszul family of ideals $\JJ$ contains $A_+$ (so that it gives a Koszul filtration)
then $A$ is Koszul.

\section{Koszul configurations of points}\label{main-sec}

Recall that we call a set $S$ of (distinct) points in $\P^n$ {\it Koszul} if
the homogeneous coordinate ring $A_S$ of $S$ is Koszul.

The next Theorem \ref{main-tech-thm} is our main technical result from which 
Theorem \ref{mainthm} will follow. 

Let $S\sub\P^n$ be a finite set of distinct points. Recall that $\LL(S)$ denotes the lattice of all subsets of $S$ obtained as intersections of $S$ with linear subspaces in $\P^n$.
For every subset $T\sub S$ let us denote by $J_T\sub A_S$ the kernel of the natural homomorphism
$A_S\to A_T$. Note that we have $A_{\emptyset}=k$, so that $J_{\emptyset}=(A_S)_+$.
  
\begin{thm}\label{main-tech-thm}  Let $S\sub\P^n$ be a finite set of distinct points. Let 
$\La\sub\LL(S)$ be a subset with the following property

\noindent
$(\ast\ast)$ for every $T\in\La$ such that $T\neq S$ there exists $T'\in\La$ such that
$T'$ is linearly independent, $T\cap T'=\emptyset$, $T\cup T'$ is an element of $\La$, and
$\dim\span(T\cup T')=\dim\span(T)+1$.

Then $\JJ=\{J_T\ |\ T\in\La\}$ is a Koszul family of ideals in $A_S$.
\end{thm}

Let us set $V=H^0(\P^n,\OO(1))$.

\begin{lem}\label{tech-lem} 
In the situation of Theorem \ref{main-tech-thm} let
us trivialize $\OO_{\P^n}(1)$ near $S$ and consider the linear map
$$e_S:V\otimes V\to k^S$$
obtained as the composition of the product map $V\otimes V\to H^0(\P^n,\OO(2))$
with the evaluation map at $S$. Assume that $\La$ contains at least one element $T$ such that
$T\neq S$. Then 

\noindent
(a) $e_S$ is surjective.

\noindent
(b) For every $T\in\La$ consider the subspace $T^{\perp}\sub V$ consisting of
linear forms vanishing at $T$. Then 
the natural map
$$e_{S,T}:V\otimes T^{\perp}\to k^{S\setminus T}$$
induced by $e_S$, is surjective.
\end{lem}

\Pf . Let us start with (b). Without loss of generality we can assume that $\span(S)=\P^n$.
We are going to use induction in codimension of $\span(T)$.
If $\codim\span(T)=0$ then $T=S$ and the assertion is trivial. Assume now that
$\codim\span(T)>0$. Then $T\neq S$, so we can apply $(\ast\ast)$ to find $T'\in\La$ disjoint from
$T$, such that $T'$ is linearly independent and $U=T\sqcup T'$ is again an element of $\La$
with smaller codimension of the span. Hence, by induction assumption
the map $e_{S,U}$ is surjective. Now from the commutative diagram
\begin{diagram}
V\otimes (T\cup T')^{\perp} & \rTo{} & k^{S\setminus (T\cup T')}\\
\dTo{} & & \dTo{}\\
V\otimes T^{\perp} & \rTo{} & k^{S\setminus T}
\end{diagram}
we conclude that it suffices to prove surjectivity of the map
$$V\otimes T^{\perp}\to k^{T'}.$$
Let $v\in T^{\perp}$ be an element that does not vanish at any point of $T'$ (such $v$ exists
since $T\cap T'=\emptyset$). Then $V\otimes (k\cdot v)\to k^{T'}$ is surjective since
$T'$ is linearly independent. Therefore, the above map is also surjective.

Now let us prove (a). Condition $(\ast\ast)$ implies that there exists a linearly independent
element $T\in\La$. Since by (a) the map $e_{S,T}$ is surjective, the commutative diagram
\begin{diagram}
V\otimes T^{\perp} & \rTo{} & k^{S\setminus T}\\
\dTo{} & & \dTo{}\\
V\otimes V & \rTo{} & k^S
\end{diagram}
reduces surjectivity of $e_S$ to that of the map
$$V\otimes V\to k^T.$$
But this immediately follows from the linear independence of $T$. 
\ed

\noindent
{\it Proof of Theorem \ref{main-tech-thm}}. If $\La=\emptyset$ or $\La=\{S\}$ then
the assertion is trivial, so we can assume that there exists $T\in\La$ such that $T\neq S$.
Hence, Lemma \ref{tech-lem} is applicable. Trivializing $\OO_{\P^n}(1)$ near $S$
and using part (a) of this lemma we can identify the algebra
$A_S$ with $k\oplus V\oplus k^S\oplus k^S\oplus\ldots$.
On the other hand, using part (b) we can derive that for every $T\in\La$
the ideal $J_T\sub A_S$ is generated in degree $1$. Indeed,
$(J_T)_1=T^{\perp}$ and $(J_T)_n=k^{S\setminus T}\sub k^S$ for $n\ge 2$, so
our claim follows from surjectivity of $e_{S,T}$. 
Now let $T\in\La$ be an element such that $T\neq S$
and let $T'\in\La$ be an element provided by condition $(\ast\ast)$. 
Let us choose $v\in T^{\perp}$ such that $v$ does not vanish at any point of $T'$.
Then
$$T^{\perp}=(T\cup T')^{\perp}+k\cdot v.$$
Since $J_T$ is generated in degree $1$, this implies that
$J_T=J_{T\cup T'}+A_S v$. 
Thus, the map $A_S(-1)\to J_T/J_{T\cup T'}: f\mapsto f\cdot v$
is surjective. The kernel of this map is the annihilator of 
$v\in A_S/J_{T\cup T'}$ which is easily seen to be $J_{T'}$. Hence, we get an
isomorphism
$$J_T/J_{T\cup T'}\simeq A_S/J_{T'}(-1).$$
\ed

\noindent
{\it Proof of Theorem \ref{mainthm}}.
First, let us show that every subset 
$\La\sub\LL(S)$, maximal among subsets satisfying $(\ast)$, also satisfies $(\ast\ast)$. 
Let $\La$ be such a maximal subset and let $T\in\La$ be any element with $T\neq S$. Then there exists 
a nonempty linearly independent element $T'\in\La$ such that $T\cup T'\in\La$ and $T\cap T'=\emptyset$.
Take a point $p\in T'$ and consider $T_1\in\LL(S)$ defined by
$$T_1=S\cap\span(T\cup\{p\}).$$
Note that $T_1\sub T\cap T'$ and $\dim\span T_1=\dim\span T+1$.
Also, since $T'$ is linearly independent, the sublattice $\LL(T')\sub\LL(S)$ consists of all subsets of $T'$.
Let us set 
$$\wt{\La}=\La\cup\LL(T')\cup\{T_1\}\sub\LL(S).$$ 
It is easy to see that condition $(\ast)$ still holds for $\wt{\La}$.
Hence, we should have $\wt{\La}=\La$. Now set $T''=T_1\setminus T$. Then $T''\sub T'$, so $T''\in\La$ and
$T''$ is linearly independent.
Thus, we have $T\cap T''=\emptyset$, $T\cup T''=T_1\in\La$ and $\dim\span(T\cup T'')=\dim\span(T)+1$,
as required in $(\ast\ast)$. 

Thus, we can assume that $\La$ satisfies $(\ast\ast)$. Note also that if $|\La|>1$
then $\La$ necessarily contains a linearly independent element $T$. Enlarging $\La$ if needed we can assume
that $\La$ contains the entire sublattice $\LL(T)\sub\LL(S)$, so that $\emptyset\in\La$.
Now Theorem \ref{main-tech-thm} implies that the family of ideals $\{ J_T\ |\ T\in\La'\}$ gives a Koszul filtration
of $A_S$ (recall that $J_{\emptyset}=(A_S)_+$).
Hence, the algebra $A_S$ is Koszul (see Corollary \ref{Koszul-fam-cor}).
\ed

\noindent{\bf Final remarks.}

\noindent
{\bf 1.} Koszulness of $S$ in Theorem \ref{mainthm} implies the same property for any subset $S'\sub S$ by Corollary \ref{descent-cor}, since $S$ imposes independent conditions on quadrics (see Lemma \ref{tech-lem}(a)). One can also check directly that the assumptions of this theorem are satisfied for any
$S'\sub S$ (where $S'\neq\emptyset$). Namely, let $\La'\sub\LL(S')$ be the image of $\La$ under the
natural intersection map $\LL(S)\to \LL(S')$. It is easy to see that $|\La'|>1$.
We claim that condition $(\ast)$ holds for $\La'$. Indeed, for any $t\in\La'$ such that $t\neq S'$
consider an element $T\in\La$ such that $T\cap S'=t$ and $T$ is maximal with respect to inclusion
among all such elements of $\La$. Applying condition $(\ast)$ to $T$ we get a linearly independent subset $T'\in\La$ such that $T\sqcup T'\in\La$. Now we set $t'=T'\cap S'$. It is clear that $t'$ is
linearly independent and that $t\sqcup t'\in\La'$. Also, $t'\neq\emptyset$ since otherwise $T$ would not be maximal. Hence, condition $(\ast)$ holds for $\La'$.

\noindent {\bf 2.} Every maximal subset $\La\sub\LL(S)$ satisfying $(\ast)$ is stable 
under taking intersections. Indeed, first, we notice that $\La$ contains all subsets of
its linearly independent elements. Now assume that $T_1, T_2\in \La$ but 
$T_1\cap T_2\not\in\La$. We can assume also that for every $\wt{T}_1\in\La$ such that
$T_1\sub\wt{T}_1$ and $T_1\neq\wt{T}_1$ one has $\wt{T}_1\cap T_2\in\La$. We
are going to lead this to contradiction by showing that $\La\cup\{T_1\cap T_2\}$ still
satisfies $(\ast)$. Indeed, let us choose a linearly independent element $T'_1\in\La$ such that $T'_1\neq\emptyset$ and $T_1\sqcup T'_1\in\La$. Then 
$$(T_1\sqcup T'_1)\cap T_2=(T_1\cap T_2)\sqcup (T'_1\cap T_2)\in\La.$$
Since $T_1\cap T_2\not\in\La$, this implies that $T'_1\cap T_2\neq\emptyset$.
Hence, condition $(\ast)$ holds for $T_1\cap T_2$.

\noindent {\bf 3.}
Let $S\sub\P^n$ be a set of $s$ points imposing independent conditions on quadrics. Then $S$ can be Koszul only when $s\le 1+n+n^2/4$. Indeed, the Hilbert series of the quotient of
$A_S$ by a generic element of degree $1$ will be $h(z)=1+nz+(s-n-1)z^2$. Since such a quotient
is still a Koszul algebra the series $h(-z)^{-1}$ should have positive coefficients. This easily implies
the required inequality (see the proof of Thm 3.1 of \cite{CTV}). Conversely, 
for every $s\le 1+n+n^2/4$ our Theorem \ref{mainthm}
gives examples of Koszul sets of $s$ points in $\P^n$ imposing independent conditions on quadrics. 
Indeed, for $m\le n+1$ consider 
$m$ linearly independent hyperplanes $H_1,\ldots,H_m$ and let $S_i$ be a linearly independent set 
spanning  $\cap_{j\neq i}H_j$, so that $|S_i|=n-m+2$. Theorem \ref{mainthm} implies
that if $S_i$ does not intersect $H_i$ for all $i$ then $S=\cup_i S_i$ is Koszul (one should take
as $\La$ all subsets of the form $S_{i_1}\cup\ldots\cup S_{i_r}$).
The number of elements in this set is $s=m(n-m+2)$. If $n$ is even then we get the
maximum of $s$ for $m=1+n/2$ which gives $s=(1+n/2)^2=1+n+n^2/4$. If $n$ is odd then for 
$m=(n+1)/2$ we get
$s=(n+1)(n+3)/4=3/4+n+n^2/4$ which is equal to $1+n+[n^2/4]$.
Passing to subsets we get required examples of Koszul sets imposing independent conditions
on quadrics for arbitrary $s\le 1+n+n^2/4$.
These examples can be used to give a different proof of Thm. 4.1 of \cite{CTV}.

\noindent {\bf 4.}
The above configurations are rather special. On the other extreme Corollary \ref{maincor}
describes the most generic type of configurations whose Koszulness follows from Theorem 
\ref{mainthm}. The following corollary describes a wider class of Koszul configurations.

\begin{cor}\label{lessgen-cor} Assume that $S=S_1\sqcup S_2\sqcup\ldots\sqcup S_m$, where every $S_i$ is linearly independent and for every $i\le j$ one has
$$S_i\sqcup S_{i+1}\sqcup\ldots\sqcup S_j\in\LL(S).$$
Then $S$ is Koszul.
\end{cor}

Indeed, in this case we can take $\La=\{S_i\sqcup S_{i+1}\sqcup\ldots\sqcup S_j\ |\ 1\le i\le j\le m
\}$.
For example, let $H, H'$ be a pair of hyperplanes in $\P^n$ and let $S_1\sub H\setminus H'$,
$S_3\sub H'\setminus H$ be a pair of linearly independent subsets of $\le n-1$ points.
Then for any linearly independent subset $S_2\sub H\cap H'\setminus (\span(S_1)\cup \span(S_3))$
the set $S=S_1\sqcup S_2\sqcup S_3$ 
is Koszul. Note that the number of elements in such a configuration is $\le 3(n-1)$. For example,
for $n=4$ in this way we get a $32$-dimensional family of Koszul configurations
of $9$ points in $\P^4$.

Here is another class of Koszul configurations.

\begin{cor}\label{D4-cor} Assume that $S=S_0\sqcup S_1\sqcup S_2\sqcup S_3$, where
every $S_i$ is linearly independent. Assume also that $S_i\in\LL(S)$ for all $i$ and
$$S_0\sqcup S_i\sqcup S_j\in \LL(S) \text{ for }1\le i<j\le 3.$$
Then $S$ is Koszul.
\end{cor}

In this case we can take $\La$ consisting of the subsets $S_i$, $0\le i\le 3$,
$S_0\sqcup S_i$, $1\le i\le 3$, $S_0\sqcup S_i\sqcup S_j$, $1\le i<j\le 3$, and $S$.
For example, one can take three linearly independent hyperplanes $H_1,H_2,H_3\sub\P^n$
and for each $1\le i\le 3$ take $S_i\sub H_i$ to be a linearly independent subset of $\le n-2$ points,
such that $S_i$ does not intersect two other hypeplanes. Then for any linearly independent
subset $S_0\sub H_1\cap H_2\cap H_3\setminus\cup_{i=1}^3\span(S_i)$ the set
$S=\cup_{i=0}^3 S_i$ is Koszul. The number of elements in such a configuration is $\le 4(n-2)$.

\noindent {\bf 5.}
Theorem \ref{mainthm} does not describe all Koszul configurations of points imposing 
independent conditions on quadrics. Indeed, it is easy to see that it is applicable
to a set of $s$ points in general linear position in $\P^n$ only when $s\le 2n$.
However, there exist Koszul sets of $s$ points in general linear position in $\P^n$
and imposing independent conditions on quadrics
for every $s\le 1+n+n^2/4$ (see Thm. 4.2 of \cite{CTV}).

\noindent{\bf 6.}
In conclusion we observe that
Theorem \ref{mainthm} could be useful in attacking the following

\vspace{2mm}

{\bf Conjecture.}
{\it A set of distinct points in $\P^n$ is Koszul provided it intersects any linear subspace
of dimension $r$ in $\le 2r$ points.}

\vspace{2mm}

It is well known that the homogeneous coordinate algebra of such a set is quadratic
(by the result of \cite{L} or by Thm. 2.1 of \cite{G}). Moreover, the set of points with these properties
imposes independent condtions on quadrics. Hence, by Corollary \ref{descent-cor}
it suffices to consider the case of $2n$ points in $\P^n$.
Most of configurations considered in the above conjecture
admit a partition into two subsets as in Corollary \ref{maincor}. For example, this can be easily
checked for $n\le 3$. The first case
when this is not true is for $8$ points in $\P^4$. The following counterexample was found by Eric Rain.
Let $e_1,\ldots, e_5$ be a basis in $k^5$ and let
$S$ be the image in $\P^4$ of the set of nonzero vectors $\{e_1,\ldots,e_5,x,y,z\}$, where 
$$x=e_1-e_3+e_4,\ y=e_1+e_2+e_4-e_5,\ z=e_1+e_2+e_3.$$
Then $S$ cannot be partitioned into two quadruples of points as in Corollary \ref{maincor}.
However, one can check that for $8$ points in $\P^4$ the above conjecture follows from Corollary
\ref{lessgen-cor} with $m\le 3$. For example, for the set $S$ above one can take $S_1=\{e_2,e_3,e_4\}$,
$S_2=\{e_1,x,z\}$, $S_3=\{e_5,y\}$.
It would be interesting to investigate the situation for $n\ge 5$.




\begin{thebibliography}{99}
\bibitem{BGS} A.~Beilinson, V.~Ginzburg, W.~Soergel,
{\it   Koszul duality patterns in representation theory},
Journ.\ Amer.\ Math.\ Soc. 9 (1996) 473--527.
\bibitem{BF}  J.~Backelin, R.~Fr\"{o}berg,
{\it  Koszul algebras, Veronese subrings and and rings with linear
  resolutions}, Revue Roumaine Math.\ Pures Appl.  30 (1985), 85--97.
\bibitem{CRV} A.~Conca, M.~E.~Rossi, G.~Valla,
{\it Gr\"obner flags and Gorenstein algebras},
Compositio Math. 129 (2001), 95--121.
\bibitem{CTV} A.~Conca, N.~V.~Trung, G.~Valla,
{\it Koszul property for points in projective spaces},
Math.\ Scand. 89 (2001), 201--216.
\bibitem{FV} M.~Finkelberg, A.~Vishik.
{\it The canonical ring of generic curve of genus $\ge 5$ is Koszul},
J.\ Algebra 162 (1993), 535--539
\bibitem{G} M.~Green, {\it The Eisenbud-Koh-Stillman conjecture on linear syzygies}, Invent. Math. 136
(1999), 411--418.
\bibitem{K} G.~Kempf, {\it  Syzygies for points in projective space},
Journ.\ Algebra 145 (1992), 219--223.
\bibitem{L} Ngau Lam, {\it On the conjecture of the resolution of a finite set of points in projective space},
First International Tainan-Moscow Algebra Workshop (Tainan, 1994), 229--235.
\bibitem{Piont} D.~Piontkovski, {\it  Noncommutative Koszul filtrations}, preprint math.RA/0301233.
\bibitem{P} A.~Polishchuk.
 {\it  On the Koszul property of the homogeneous coordinate ring
   of a curve},  J. Algebra 178 (1995), 122--135.
\bibitem{PP} A.~Polishchuk, L.~Positselskii, {\it Quadratic algebras}, preprint, 1996.
\bibitem{Pos} L.~Positselskii, {\it Koszul property and Bogomolov's conjecture}, Harvard Ph.~D.~thesis, 1998, available at http://www.math.uiuc.edu/K-theory/0296.
\bibitem{Pr}  S.~Priddy,  {\it Koszul resolutions},
Trans. AMS 152 (1970), 39--60.
\bibitem{ST} B.~Shelton, C.~Tingey, {\it On Koszul algebras and a new construction of
Artin-Schelter regular algebras}, J. Algebra 241 (2001), 789--798.

\end{thebibliography}
\end{document}